\documentclass[11pt]{amsart}
\usepackage{amssymb}
\usepackage{url}
\usepackage{hyperref}

\newcommand{\F}{\mathbb{F}}
\newcommand{\Z}{\mathbb{Z}}
\newcommand{\Q}{\mathbb{Q}}

\newcommand{\GL}{\mathrm{GL}}
\newcommand{\SL}{\mathrm{SL}}
\newcommand{\Sp}{\mathrm{Sp}}
\newcommand{\G}{\Gamma_n}

\newcommand{\calg}{{\mathcal G}}

\newtheorem{theorem}{Theorem}[section]

\theoremstyle{remark}

\theoremstyle{example}

\def\gp#1{\langle \hspace*{.2mm} #1 \hspace*{.25mm} \rangle}

\newcommand{\cha}{\mathrm{char} \hspace{.5pt}}

\newcommand{\abk}{\allowbreak}

\begin{document}

\title[Practical computation with linear groups
over infinite domains]{Practical computation with linear groups
over infinite domains}

\author{A.~S.~Detinko}
\address{School of Computer Science,
University of St Andrews\\
North Haugh\\
\newline St Andrews KY16 9SX\\
UK}
\email{ad271@st-andrews.ac.uk}

\author{D.~L.~Flannery}
\address{School of Mathematics, Statistics 
 and Applied Mathematics\\
National \newline University of Ireland, Galway\\
University Rd\\
Galway H91TK33\\
Ireland}
\email{dane.flannery@nuigalway.ie}

\begin{abstract}
We survey recent progress in computing with finitely 
generated linear groups over infinite fields, describing 
the mathematical background of a methodology applied to design 
practical algorithms for these groups. Implementations of 
the algorithms have been used to perform extensive computer experiments.
\end{abstract}

\maketitle

\section{Introduction}

\subsection{Motivation}
Linear groups (synonymously, matrix groups) have been studied from the 
beginning of group theory. Matrices afford a convenient
representation of groups that frequently arise in algebra, geometry,
number theory, topology, and theoretical physics. 
Enhancements of technology and computer algebra systems 
have initiated a new phase in this classical subject, concerned with the
design and implementation of algorithms for practical computation.

Computing with matrix groups over finite fields is
well-established~\cite{OBriensurvey}. 
The situation for linear groups over infinite domains is 
less advanced. Consequently we are motivated to obtain efficient methods, 
algorithms, and software 
for computing in this class of groups.  

\subsection{Representing linear groups in a computer}

Input to any algorithm should be a finite set. Thus, in the first instance,
we consider finitely generated linear groups. 
Certain linear groups that are not finitely generated can 
still be designated by a finite set---say, of polynomials, 
in the case of linear algebraic groups. Whereas an arbitrary linear group need 
not be finitely generated or algebraic, these are two major classes covering
many applications. 

Finitely generated linear groups are amenable to symbolic computation. 
Let $\F$ be a field of characteristic $p \geq 0$, and suppose that 
$G = \gp{S}$ where $S = \{g_1, \ldots , g_r\} \subseteq  \GL(n, \F)$. 
Then $G$ is defined over a finitely generated extension of the prime subfield of 
$\F$. The classification of such field extensions implies that $G$ is a subgroup of 
$\GL(n, \mathbb{L})$, where $\mathbb{L}$ is a finite degree extension 
of $\mathbb{P}({\rm x}_1, \ldots , {\rm x}_m)$, $\mathbb{P}$ is a number field 
or finite field $\F_q$ of size $q$ for some $p$-power $q$, and the ${\rm x}_i$
are algebraically independent indeterminates. This means that essentially 
we only have to deal with the aforementioned categories of fields. 
All of these are supported by the computer algebra system {\sc Magma}~\cite{Magma}.  

We could restrict the ground domain to the subring $R \subseteq \F$ 
generated by all entries of the $g_i$ and $g_i^{-1}$. After replacing the original
field by such a ring, we apply congruence homomorphism techniques to 
transfer computing over $R$ to computing over a quotient ring $R/\rho$.
If $\rho$ is a maximal ideal then $R/\rho$ is a finite field, and in  
that event the computational complexity is ameliorated by avoiding work over
an infinite ring. We
also gain access to the machinery for matrix groups 
over finite fields. See \cite[Section~2]{Recognition} for details.

\subsection{Properties of linear groups}

We rely on classical theory of linear groups~\cite{DixonBook,SuprunenkoI,WBook}. 
Two basic properties are crucial in our endeavours. One of these  provides 
background for the computational methods; the other steers our overall strategy. 

First, we recall that each finitely generated linear group $G$ is residually 
finite. Moreover, $G$ is `approximated' by matrix groups of the same degree over 
finite fields. This approximation is effected by congruence homomorphisms 
$\varphi_{\rho} : \GL(n, R) \rightarrow \GL(n, R/\rho)$. 
Since each non-zero element of $R$ is absent from at least one ideal,
and $R/\rho$ is a finite field if we choose $\rho$ to be maximal, 
the congruence images $\varphi_{\rho}(G)$ realize the finite approximation. 

A famous result of J.~Tits~\cite{Tits72} asserts that each finitely generated 
linear group over a field either is solvable-by-finite (virtually solvable), 
or contains a free non-abelian subgroup. The Tits alternative thereby divides 
finitely generated linear groups into two very different classes which
require separate treatment.

\section{Computing with virtually solvable groups}\label{ComputingWithSFGroups}

\subsection{Method of finite approximation} 

Our techniques for computing with solvable-by-finite 
groups are broad-based and uniform, enabling us to solve 
a range of problems by similar algorithms. Underlying 
these features are deep results 
about the congruence subgroup
$G_{\rho}:=G\cap \ker \, \varphi_\rho$.
\begin{theorem}\label{DefineWHom}
There exist maximal ideals 
$\rho$ of $R$ such that
\begin{itemize}
\item[{\rm (i)}] All torsion elements of $G_{\rho}$ are unipotent. 
In particular, $G_{\rho}$ is torsion-free if $\cha R = 0$.
\item[{\rm (ii)}] If $G$ is solvable-by-finite then $G_{\rho}$ is 
unipotent-by-abelian as long as one of the following holds:
$\cha R >n$; 
$\cha R =0$ and $\cha (R/\rho) >n$;
$R$ is a Dedekind domain of characteristic zero and
$p\in \rho\setminus\rho^{p-1}$ for some odd prime $p$.   
\end{itemize}
\end{theorem}
See \cite[Chapter 4]{WBook} or \cite[Section~2]{Recognition}
for a proof of Theorem~\ref{DefineWHom}~(i).
Proofs, and extra conditions on $R$ and $\rho$ 
guaranteeing the outcome in 
Theorem~\ref{DefineWHom}~(ii),
are given in \cite{WJAlgebra}.

Our method begins by selecting $\rho$ according to the strictures 
of Theorem~\ref{DefineWHom}, and 
computing the congruence image $\varphi_{\rho}(G) \leq \GL(n, R/\rho)$.
Then we examine the structure of $G_{\rho}$. 
We call $\varphi_{\rho}$ for $\rho$ as in Theorem~\ref{DefineWHom} 
a \emph{W-homomorphism}. 
Algorithms to compute W-homomorphisms $\varphi_\rho$
and their corresponding 
congruence images $\varphi_{\rho}(G)$ were developed in \cite{Tits,Recognition}. 
These compute normal 
generators of $G_{\rho}$, i.e., a finite set $N \subseteq G$ such that 
$G_{\rho} = \gp{N}^G$. The set $N$ is found by means of a 
presentation of $\varphi_{\rho}(G)$, computed using 
algorithms for matrix groups over finite fields~\cite{Baaetal,OBriensurvey}. 
For our purposes, any relevant information about $G_\rho$ can be deduced
from $N$; 
the full normal closure $\gp{N}^G$ is not needed.

\subsection{Recognizing the type of a matrix group}

Armed with practical methods, we proceed to the 
development of algorithms. 
Given $S\subseteq \GL(n,\F)$ 
we must first recognize the `type' of $G = \langle S \rangle$. 
Once this is done, $G$ can be investigated using tools that are most 
appropriate for the group type. Below we note algorithms to recognize 
the type of $G$ (each of which requires selection of a single 
 W-homomorphism $\varphi_{\rho}$). 
These algorithms additionally justify that the relevant  
problems are decidable for finitely generated linear groups over 
infinite fields.

\subsubsection{Finiteness} 
In characteristic zero, $G$ is finite if and only if
$G_{\rho} = \langle N \rangle^G=1$.
If $\cha \, \F = p > 0$ then finiteness testing turns on
whether $G_{\rho}$ is a $p$-group, i.e., unipotent. See 
\cite[Section~4]{Recognition}.

\subsubsection{Virtual solvability and other properties} 
\label{TestingVirtualSolvability}
We can recognize whether $G$ is solvable-by-finite: a computational 
realization of the Tits alternative. For this it is enough to test 
whether $G_{\rho}$ is unipotent-by-abelian, i.e., conjugate to 
a block-triangular group with all main diagonal blocks abelian. 
This test is carried out using manipulations with the enveloping 
algebra of $G_{\rho}$ over $\F$, as explained in \cite[Section~3]{Tits}. 
Although it decides whether a finitely generated linear group 
contains a free non-abelian subgroup, our algorithm does not construct 
one. 

Algorithms to test whether $G$ is solvable, 
(virtually) nilpotent, abelian-by-finite, or central-by-finite. 
use a mix of ideas similar to the above~\cite[Section~5]{Tits}.

\subsection{Investigating the structure of linear groups}

\subsubsection{Finite groups} 
If $G$ is found to be finite then we can obtain an isomorphic copy 
 over a finite field $\F_q$. In characteristic zero, 
$G \cong \varphi_{\rho}(G)$ 
for any W-homomorphism $\varphi_{\rho}$; in positive characteristic, 
repeated selection of $\rho$ may be needed to get an isomorphism 
$\varphi_\rho$~\cite[Section~4.3]{Recognition}. 
Algorithms for matrix groups over finite fields may then be
 applied to
$\varphi_{\rho}(G)\leq \GL(n,q)$ to answer questions about the original
group $G$. 

\subsubsection{Solvable groups}

Linear groups play a central role in the theory of infinite 
solvable groups. However, in designing algorithms for solvable linear 
groups we encounter serious obstacles, such
as lack of decidability of various problems~\cite[Chapter 9]{LennoxRobinson}.  
To further illustrate this point, we make a comparison with polycyclic groups. 
Virtually polycyclic groups are finitely generated and $\Z$-linear. 
On the other hand, finitely generated (virtually) solvable linear groups 
need not be finitely presentable, they might have subgroups that are not 
finitely generated, and they do not satisfy the maximal condition on subgroups. 
Computing becomes feasible with groups of finite Pr\"ufer rank, which
are solvable-by-finite and $\Q$-linear. Hence, we can test whether a 
finitely generated linear group $G$ over a number field $\F$ has finite
rank. Furthermore, if $G$ is (virtually) solvable then we can: 
compute the torsion-free rank (Hirsch number) of $G$, and bounds on its Pr\"ufer rank;  
test whether $|G : H|$ is finite, for a finitely generated subgroup $H$ of $G$; 
construct a generating set of the completely reducible part of $G$ (this 
includes testing whether $G$ is completely reducible or unipotent).
More generally, these algorithms work for solvable-by-finite groups $G$ over 
any field, albeit with qualifications on $G$ in positive characteristic.   
The papers \cite{Tits,FiniteRank} contain lengthier discussion of the above.

Nilpotent-by-finite linear groups are more tractable. 
Algorithms for these are given in \cite{Nilpotent}
 and \cite[Section~5]{Tits}.
Computing with polycyclic linear groups is a separate topic
(see, e.g., \cite{AssmannEick1,AssmannEick2}) beyond
the remit of our survey.

\subsubsection{Implementation} 
Many of our algorithms for virtually solvable groups were developed jointly 
with Eamonn O'Brien. 
Implementations are available in {\sc Magma}; see \cite{Infinite}. Experimental
results are reported in
\cite[Section~6]{Tits}, 
\cite[Section~4.5]{FiniteRank}, and
\cite[Section~5]{Recognition}. 

\section{Dense and arithmetic groups}\label{DenseAndArithmetic}

The methods of Section~\ref{ComputingWithSFGroups} could be developed further. 
However, to move beyond virtually solvable groups, new 
ideas are required.  

Each linear group $H$ is contained in an algebraic group, with the 
Zariski closure of $H$ being the `smallest' such overgroup. 
We will suppose that $H$ is a dense (in the Zariski topology) 
subgroup of an algebraic group. 
Note that an algorithm to compute the Zariski closure 
of a finitely generated linear group is given in \cite{Derksen}.

The most interesting case is $\Q$-groups 
$\calg \leq \GL(n, \mathbb{C})$, i.e., $\calg$ is 
defined by a set of polynomials with 
coefficients in $\Q$. For a subring $R \subseteq \mathbb{C}$, denote 
$\calg \cap \GL(n, R)$ by $\calg(R)$. Recall that 
$H \leq \calg(\Q)$ is \emph{arithmetic} if 
$H \cap \calg(\Z)$ has finite index in $H$ and in
$\calg(\Z)$. In particular, finite index subgroups of $\calg(\Z)$ are  
arithmetic. Arithmetic groups are finitely generated and dense. 
If $H\leq \calg(\Z)$ is dense but not arithmetic, then we 
call $H$ a \emph{thin matrix group} (after \cite{Sarnak}).
A major open problem is testing whether finitely 
generated subgroups of $\calg(\Z)$ are arithmetic. In \cite{DeGrDetF} 
we provide an algorithm (implemented in {\sc Magma})
to test arithmeticity when $\calg$ is solvable; showing 
that the problem is decidable with this proviso. 
The algorithm computes a generating 
set of an arithmetic subgroup in $\calg(\Z)$, compares 
its Hirsch number with that of the input $H\leq \calg(\Q)$,
and tests integrality of $H$.

\subsection{Density and computing with linear groups}

Most linear groups are not virtually solvable~\cite{Aoun,Epstein}.
So we cannot expect to handle every finitely generated linear group
$H$ that is not virtually solvable by a single uniform method. Selecting 
one ideal at a time might not suffice for all problems. 

We are viewing $H$ as a subgroup of some algebraic 
$\Q$-group $\calg$, which may be assumed semisimple by a 
standard reduction \cite[Chapters~3 and 4]{WillemBook}.
Since $H$ should be dense in $\calg$, density testing is a preliminary 
task. A deterministic algorithm to test density of $H$ 
is given in \cite{RivinIII}, together with a Monte-Carlo algorithm that 
tests density 
of $H \leq \abk \calg(\Z)$ for $\calg = \SL(n, \mathbb{C})$ or 
$\Sp(n, \mathbb{C})$; see also \cite[Section~3.2]{Density}. 
These algorithms have been implemented in {\sf GAP}~\cite{GAP} (see \cite{OWGAPDoc}).

\subsection{From finite to strong approximation}
\label{FromFiniteToStrong}

We have expanded the congruence homomorphism methodology to cover 
dense subgroups. 
For certain $\calg$, and $H \leq \calg(\Z)$ 
dense in $\calg$, a celebrated result known as 
\emph{the strong approximation theorem}~\cite[Window~9]{LubotzkySegal} 
enables us to compute all congruence quotients of $H$
modulo primes.
Both $\SL(n, \mathbb{C})$ and $\Sp(n, \mathbb{C})$ 
are suitable examples of such $\calg$; from now on $\calg$ 
stands for either of these two groups. 
Strong approximation implies that if $H \leq \calg(\Z)$ 
is dense then 
$\varphi_p(H) = \varphi_p(\calg(\Z))$ for all but finitely
many primes $p$. Denote the set of these exceptional primes by 
$\Pi(H)$. We have developed practical algorithms to compute $\Pi(H)$, 
thus realizing strong approximation computationally; see
\cite[Section~3.2]{Density}, \cite{DensityFurther}, \cite{SATGeneral}. 
Our methods for computing $\Pi(H)$ draw on classifications 
of maximal subgroups in $\SL(n, p)$ and $\Sp(n, p)$, and subgroups of $\GL(n, p)$ with a known
transvection. Actually, once we have $\Pi(H)$ we can find \emph{all} congruence quotients of $H$~\cite{DensityFurther,SATGeneral}. 

\subsection{From density to arithmeticity}

Let $n > 2$ and $H\leq \calg(\Z)$ be dense. Then $H$ lies in a 
unique `minimal' arithmetic group $\mathrm{cl}(H)$, 
namely the intersection of all arithmetic groups in $\calg(\Z)$ 
containing $H$. Algorithms for arithmetic subgroups of 
$\calg(\Z)$ can therefore be used to study dense subgroups as well. 

We gain much mileage from the fact that 
$\G := \calg(\Z)$ has the \emph{congruence subgroup property}: each 
arithmetic group $H$ in $\G$ contains a principal congruence 
subgroup (PCS), which is the kernel of a congruence homomorphism 
$\varphi_m : \G \rightarrow \GL(n, \Z_m)$
for some $m$. Here $\Z_m = \Z/m\Z$, and $m$ is called the 
\emph{level} of the PCS. 
The \emph{maximal} PCS of $H$ is unique, and its level 
$M = M(H)$ is defined to be the level of $H$. Similarly, for dense
$H \leq \G$, we assign $M(H)$ as 
the level of $\mathrm{cl}(H)$. 

\subsection{Computing via the congruence subgroup property}
\label{ComputingViaCSP}

The bedrock of
our method for computing with dense groups is the 
congruence subgroup property. It splits our method into two overlapping
parts: finding $M(H)$, and computing in $\GL(n, \Z_m)$. 

\subsubsection{Computing the level} 
 
The set $\pi(M)$ of prime divisors of $M(H)$ coincides with $\Pi(H)$, 
besides minor exceptions for $n = 3$, $4$ and $p = 2$ 
(which are dealt with separately); see \cite[Section~2.4]{Density}. 
Thus, the strong approximation algorithms cited in 
Section~\ref{FromFiniteToStrong} 
yield $\pi(M)$. We can also compute the largest power of $p$ 
dividing $M(H)$ for each $p\in \pi(M)$.
These two steps constitute the procedure ${\tt LevelMaxPCS}$, which
accepts $\Pi(H)$ and a generating set $S$ of a dense group $H\leq \G$,
and returns its level. 

\subsubsection{Computing with matrix groups over $\Z_m$} 
Algorithms for subgroups of $\GL(n, \Z_m)$ have intrinsic value. 
We reduce computing to the situations of 
matrix groups over finite fields, and groups of prime-power order. 
Two major steps in the reduction are as follows. 
Say $m = p_1^{k_1} \ldots p_t^{k_t}$ where the 
$p_i$ are distinct primes and all $k_i$ are non-zero. 
Then (essentially by the Chinese Remainder Theorem)
\begin{itemize}
\item[(i)] $\GL(n, \Z_m) \cong \GL(n, \Z_{p_1^{k_1}}) 
\times \cdots \times \GL(n, \Z_{p_t^{k_t}})$
\item[(ii)] 
$\GL(n, \Z_{p^k}) / K \cong \GL(n, p)$, 
where $K = \{ h \in \GL(n, \Z_{p^k}) \; | \; 
h \equiv 1_n \allowbreak \hspace{3pt} 
\mathrm{mod} \hspace{3pt} p^{k - 1}\}$ is a $p$-group. 
\end{itemize}
We also use the fact $K\cap G$ almost always does not have a 
proper supplement in $G$, for $G=\SL(n, \Z_{p^k})$ or 
$\Sp(n, \Z_{p^k})$~\cite[Theorem 2.5]{Density}.

\subsection{Algorithms for arithmetic subgroups}

Let $H\leq \G$ be arithmetic.
In the application of Section~\ref{ComputingViaCSP} to 
designing algorithms for $H$, the main steps are 
${\tt LevelMaxPCS}$,
and computing with matrix groups over finite rings $\Z_m$. 
One example is the membership test 
${\tt IsIn}(g, H)$ which determines whether $g \in \G$ 
is in $H$; it merely checks whether 
$\varphi_{M}(g)\in \varphi_{M}(H)$.
We emphasize that our results imply decidability of 
membership testing in arithmetic groups in $\G$. 
An associated algorithm computes $|\G : H|$. Although the index 
could be calculated in the congruence image, i.e.,
as $|\varphi_{M}(\G) : \varphi_{M}(H)|$, in practice  
$|\G : H|$ is found as a byproduct of 
computing $M$~\cite[Section~2.4.2]{Density} (see 
\cite[Section 6]{Arithm} and \cite[Section 4]{Density}). 
Since membership testing and computing the index are both 
decidable, an arithmetic group $H\leq \G$ is 
`explicitly given' as per \cite{GSI}. 
Other notable algorithmic problems for arithmetic 
subgroups are therefore decidable too.

\subsection{Further computation with arithmetic subgroups}

\subsubsection{Structural analysis}

Arithmetic groups are matrix groups defined over rings, and so 
their (sub)normal structure is of interest. 
The procedure 
${\tt IsSubnormal}(H)$ tests whether $H$ is subnormal in 
$\G$; ${\tt Normalizer}(H)$ computes a generating set 
of the normalizer of $H$ in $\G$; 
${\tt NormalClosure}(B)$ computes a generating set of the normal 
closure in $\G$ of the group generated by $B\subset \G$. 
Other procedures are given in \cite[Section~3.2]{Arithm}. Many
more algorithms could be developed along these lines.

\subsubsection{The orbit-stabilizer problem} 
Let $n>2$ and $H\leq \SL(n, \Z)$ be arithmetic.
Given $u, v \in \Q^n$, the procedure ${\tt Orbit}(u, v)$ tests whether there is
$g \in \SL(n, \Z)$ such that $gu = v$, and computes $g$ if it exists. 
${\tt Stabilizer}(H, u)$ returns the (finitely generated)
stabilizer of $u$ in $H$. 
Both procedures solve the related orbit and stabilizer problems for the congruence
image over $\Z_{M}$ 
and for the maximal PCS in $H$ acting on $\Q^n$. The
outputs are then combined. See \cite[Section 4]{Arithm}.

\subsection{Experiments}
\label{ExperimentsAndApplications}

The algorithms of this section are joint work with Alex-ander 
Hulpke.  Below we review some experiments illustrating our
{\sf GAP} implementation of the algorithms 
and their practicality; 
see \cite{Density, DensityFurther, SATGeneral} for more. 

\subsubsection{}  
Integral representations of the fundamental group 
$\gp{x, y, z \, | \abk \, zxz^{-1} = xy, zyz^{-1} =yxy}$
of the figure-eight knot complement are constructed in
\cite{LongReidI}. 
For non-zero $T \in \Z$, let $\beta_T(x) = X_T$ and 
$\beta_T(y) = Y_T$ where
\[
X_T = {\small \left[\begin{array}{crc} -1+T^3 & -T & T^2
\\
0 & -1 & 2T
\\
-T & \phantom{-}0 & 1 \end{array} \right]}, \quad
Y_T = {\small \left[\begin{array}{ccr} -1 & 0 & 0
\\
-T^2 & 1 & -T
\\
\phantom{-}T & 0 & -1 \end{array} \right]}.
\]
Then  $\beta_T$ is a homomorphism 
and $\beta_T(\gp{x,y})\leq \SL(3,\Z)$ is arithmetic.
Construction of these representations was 
motivated by long-standing problems; such as 
the conjecture that each arithmetic group in $\SL(n, \Z)$ 
has a $2$-generator finite index subgroup.
The conjecture has been settled affirmatively~\cite{Meiri}.
Still, the subgroups $\langle X_T, Y_T\rangle$ merit 
closer scrutiny. Earlier attempts to compute 
$|\SL(3, \Z) : \langle X_T, Y_T\rangle|$ 
were stymied  by the fact that this index may
be arbitrarily large.  We 
were able to compute indices using our algorithms
(see \cite[Section 4.1]{Density}). For example, let $T = 100$; 
the index 
$2^{42}3^{5}5^{25}7^{4}13{\cdot}31^{2}67{\cdot}1783$ 
and level $2^{7}5^{6}29{\cdot}67{\cdot}193$ were found in
$892.6$ seconds. 

\subsubsection{}
A second family of test groups comes from applications 
in theoretical physics. 
Let $G(d, k) = \langle U, T \rangle$ where 
\[
U = {\small \left[\begin{array}{crrc} 1 & 1 & 0 & 0
\\
0 & 1 & 0 & 0
\\
d & d & 1 & 0
\\
0 & -k & -1 & 1
\end{array} \right]},
\quad  T = {\small \left[\begin{array}{cccc} 1 & 0 & 0 & 0
\\
0 & 1 & 0 & 1
\\
0 & 0 & 1 & 0
\\
0 & 0 & 0 & 1
\end{array} \right]}.
\] 
For fourteen pairs $d, k$ of integers,
$G(d, k) \leq \Sp(4, \Z)$ 
is the monodromy group of a generalized hypergeometric 
ordinary differential equation associated to Calabi-Yau threefolds. 
Seven of these groups are arithmetic, while the rest are 
thin~\cite{Singh,SinghVenky}. 
To investigate the latter, one could attempt to construct 
arithmetic groups in $\Sp(4,\Z)$ containing them~\cite{MonodromyAppendix}. 
We successfully computed 
$\mathrm{cl}(G(d, k))$ for the seven thin 
groups~\cite[Table~3]{Density}; e.g.,
it took $25$ seconds to find the level $2^{5}3^{2}$ and the 
index $2^{17}3^{6}5^{2}$ of $G(12,7)$.

\section{Where to next?}

We outline avenues for future research.

New methods and algorithms for algebraic groups and Lie algebras 
would have an impact on computing with virtually solvable groups. Despite 
significant progress (cf.~Section~\ref{ComputingWithSFGroups}), 
key algorithmic questions are still unresolved. One of these is 
membership testing. This problem is known to be decidable
for groups of finite rank. The main challenge is handling the 
unipotent radical, which is a torsion-free nilpotent group that may
not be finitely generated.
Lie algebra methods due to P.~Hall, and computing in 
ambient solvable algebraic groups, are possible approaches. These are 
similarly promising in the design of algorithms for structural analysis of 
virtually solvable linear groups. We also expect a number of new 
algorithms for computing with (virtually) nilpotent and (virtually) 
polycyclic linear groups.

Methods based on algebraic group techniques will be  
productive in applications to non-virtually solvable groups 
(cf.~Section~\ref{DenseAndArithmetic}). Arithmeticity testing is 
open in general,
even for subgroups of $\SL(n, \Z)$. Indeed, it is not known whether the 
problem is decidable. Computing
generating sets and presentations of arithmetic subgroups 
are supplementary problems (cf. \cite[Chapter 6]{WillemBook}, \cite{Nebe}). 
Construction of 
free subgroups would aid in the study of matrix groups 
that are not virtually solvable; `large' free subgroups, 
i.e., those that are dense in the Zariski closure, are especially useful. 
Testing freeness of finitely generated linear groups is yet another 
priority. 

We await breakthroughs that apply computational methods 
to the solution of hard problems in group theory, 
other areas of mathematics,
and farther afield (cf.~Section~\ref{ExperimentsAndApplications}).
Here we point to computing linear representations of finitely 
presented groups: 
in contrast to the same problem for finite groups, much
remains to be done.

\subsection*{Acknowledgments}
We are indebted to our collaborators Willem de Graaf, Alexander Hulpke, and
Eamonn O'Brien.
We also thank Mathematisches Forschungsinstitut Oberwolfach,
and the International Centre for Mathematical Sciences, UK, for hosting our
visits under their `Research in Pairs' and `Research in Groups'
programmes. A.~S.~Detinko is supported by a
Marie Sk\l odowska-Curie Individual
Fellowship grant (Horizon 2020, EU Framework Programme for
Research and Innovation).

\medskip

\bibliographystyle{amsplain}
\def\cprime{$'$}
\providecommand{\bysame}{\leavevmode\hbox to3em{\hrulefill}\thinspace}
\providecommand{\MR}{\relax\ifhmode\unskip\space\fi MR }
\providecommand{\MRhref}[2]{%
  \href{http://www.ams.org/mathscinet-getitem?mr=#1}{#2}
}
\providecommand{\href}[2]{#2}

\bigskip

\end{document}